\documentclass[12pt]{amsart}

\usepackage[margin=1in]{geometry}

\usepackage{amsmath,amscd,amssymb,amsfonts,latexsym,wasysym, mathrsfs, mathtools,hhline,color}
\usepackage[all, cmtip]{xy}
\usepackage{csquotes}
\usepackage{url}
\usepackage{comment}

\definecolor{hot}{RGB}{65,105,225}

\usepackage[pagebackref=true,colorlinks=true, linkcolor=hot ,  citecolor=hot, urlcolor=hot]{hyperref}%make all the references and links clickable
\usepackage{ textcomp }
\usepackage{ tipa }
\usepackage{graphicx,enumerate}

\usepackage{geometry}
\theoremstyle{plain}
\newtheorem{theorem}{Theorem}[section]
\newtheorem{prop}[theorem]{Proposition}

\newtheorem{lm}[theorem]{Lemma}

\newtheorem{cor}[theorem]{Corollary}
\newtheorem{conj}[theorem]{Conjecture}

\newtheorem{thrm}[theorem]{Theorem}

\theoremstyle{definition}

\newtheorem{defn}[theorem]{Definition}

\newtheorem{rmk}[theorem]{Remark}

\newtheorem{ex}[theorem]{Example}
\newtheorem*{ex*}{Example}

\def\be{\begin{equation}}
\def\ee{\end{equation}}

\def\bt{\begin{thrm}}
\def\et{\end{thrm}}

\def\bc{\begin{cor}}
\def\ec{\end{cor}}

\def\br{\begin{rmk}}
\def\er{\end{rmk}}

\def\bp{\begin{prop}}
\def\ep{\end{prop}}

\def\bl{\begin{lm}}
\def\el{\end{lm}}

\def\bex{\begin{ex}}
\def\eex{\end{ex}}

\def\bd{\begin{defn}}
\def\ed{\end{defn}}

\newcommand{\C}{\mathbb{C}}
\newcommand{\Z}{\mathbb{Z}}

\newcommand{\tor}{\mathrm{tor}}

\newcommand{\sA}{\mathcal{A}}

\newcommand{\sL}{\mathcal{L}}

\newcommand{\rank}{\mathrm{rank }}

\begin{document}

\title[Double coverings and local systems]{Integral homology groups of double coverings and rank one $\Z$-local systems for a minimal CW complex}


%    author one information
\author{Ye Liu}
\address{Department of Pure Mathematics, Xi'an Jiaotong-Liverpool University, 111 Ren'ai Road, Suzhou, Jiangsu 215123 China}
%\curraddr{}
\email{yeliumath@gmail.com}
%\thanks{}

%    author two information
\author{Yongqiang Liu}
\address{The Institute of Geometry and Physics, University of Science and Technology of China, 96 Jinzhai Road, Hefei, Anhui 230026 China}
%\curraddr{}
\email{liuyq@ustc.edu.cn}
%\thanks{}

%    \subjclass is required.
\subjclass[2020]{Primary 55N25, 52C35}

\date{}

%    Abstract is required.
\begin{abstract}
Given a finite CW complex $X$, a nonzero cohomology class $\omega \in H^1(X,\Z_2)$ determines a double covering $X^\omega$ and a rank one $\Z$-local system $\sL_\omega$. We investigate the relations between the homology groups $H_*(X^{\omega},\Z)$ and $H_*(X,\sL_\omega)$, when $X$ is homotopy equivalent to a minimal CW complex. In particular, this settles a conjecture recently proposed by Ishibashi, Sugawara and Yoshinaga for a hyperplane arrangement complement \cite[Conjecture 3.3]{ISY22}.
\end{abstract}

\maketitle


\section{Introduction}

Finite coverings of CW complexes are classical subjects in algebraic topology. Recently, double coverings attract a lot of attentions for the topology of hyperplane arrangement complement \cite{Yos20,Suc22,ISY22}.

A finite collection $\sA$ of hyperplanes in $\C^n$ (or $\mathbb{CP}^n$) is called a complex affine (resp. projective) hyperplane arrangement. The topology of the hyperplane arrangement complement is very interesting. For instance, Dimca and Papadima \cite{DP03} and Randell \cite{Ran02} independently showed that the complement of a hyperplane arrangement is homotopy equivalent to a minimal CW complex.   A fundamental problem in the theory of hyperplane arrangements is to decide whether various topological invariants of the complement of $\sA$ are determined by the combinatorial structure of $\sA$. It is well known that the Betti numbers and the cohomology ring of the hyperplane arrangement complements are
combinatorially determined (e.g., see \cite{OT}). However, it is still an open question whether
 the Betti numbers of  a finite abelian cover  of a hyperplane arrangement
complement are combinatorially determined. This includes the Milnor fiber of a central hyperplane arrangement. See \cite{PS17} for recent progress in this direction and also \cite{Suc01} for an overview of the theory.

Yoshinaga studied the mod-2 Betti numbers of double coverings for a hyperplane arrangement complement and showed that they are combinatorially determined \cite[Theorem 3.7]{Yos20}. 
As an application, he showed that the first integral homology group of the Milnor fiber of the icosidodecahedral arrangement has 2-torsion \cite[Theorem 1.2]{Yos20}. Ishibashi, Sugawara and Yoshinaga 
\cite{ISY22} further studied the 2-torsion part of the homology groups of the double coverings and gave a refinement of Papadima and Suciu's work \cite{PS10}. Based on computations, they proposed a conjecture \cite[Conjecture 3.3]{ISY22} regarding the first integral homology group of a double covering and the first homology group of a rank one $\Z$-local system. 
In this note, we settle this conjecture in all degrees with the more general setting: the CW complex is homotopy equivalent to a minimal CW complex. As an application, we show that the integral homology groups of double coverings of hyperplane arrangement complement are combinatorially determined under certain conditions.

%This formula is  strengthened by Suciu  \cite[Proposition 8.5]{Suc22} to full generality recently.

\bigskip

Let $X$ be a finite connected CW complex.
Fix a nonzero element $\omega \in H^1(X,\Z_2)$, where $\Z_2=\Z/2\Z$. Then  $\omega$ determines a 
surjective map $\pi_1(X) \to \Z_2\cong \{\pm 1\}$. This gives a
double covering    $ X^\omega \to X$. On the other hand, the group homomorphism $\pi_1(X)\to \{\pm 1\}=\Z^\times$ also gives a rank one $\Z$-local system which we denote by $\sL_\omega$.

What is the relation between  $H_*(X^\omega,\Z)$ and $H_*(X,\sL_\omega)$?
It is easy to see that $$H_i(X^\omega, \C)\cong H_i(X,\C)\oplus H_i(X, \sL_\omega\otimes_\Z \C).$$ So the difficult part is about the torsions of the homology groups. We give a complete answer to this question when $X$ is homotopy equivalent to a minimal CW complex.

\bd Let $X$ be a connected finite CW complex. We say that the CW-structure on $X$ is minimal if the number of $i$-cells of $X$ coincides with the (rational) Betti number $b_i(X)$, for every $i\geq 0$. Equivalently,
the boundary maps in the cellular chain complex $C_*(X, \Z)$ are all the zero maps.
\ed

\bt \label{main} Let $X$ be a connected finite CW complex, which is homotopy equivalent to a minimal CW complex. 
Fix a nonzero element $\omega \in H^1(X,\Z_2)$. Then there exists a bounded complex $(E_*,\partial_*)$ of finitely generated free abelian groups
$$\cdots \to E_{i+1} \overset{\alpha_{i}}{\to} E_i \overset{\alpha_{i-1}}{\to} E_{i-1} \to \cdots $$
such that $H_i(E_*,\alpha_*) \cong H_i(X,\sL_\omega)$, $\rank E_i=b_i(X)$ and every entry in the boundary map $\alpha_i$ is divisible by 2. Then the complex $(E_*, \tfrac{1}{2}\alpha_*)$ is well defined. Moreover, we have
\begin{equation}\label{*}
H_i(X^\omega,\Z)\cong H_i(X,\Z) \oplus H _i(E_*, \tfrac{1}{2}\alpha_*).
\end{equation}
\et

\br The minimal CW complex assumption is crucial for the above theorem. For example, consider the real projective space $\mathbb{RP}^n$ with $n>1$. Its double covering space is the sphere $S^n$. For any $0<i<n$, we have that $ H_i(S^n,\Z)=0$, yet $$H_i(\mathbb{RP}^n,\Z)=\begin{cases}
\Z_2, & \text{ for } i \text{ odd },\\
0, & \text{ for } i \text{ even }.
\end{cases}$$ %For more results in this direction, see  \cite[Section 8]{Suc22}.
\er

\bc  \label{cor1} With the same assumptions and notations as in Theorem \ref{main}, we have that
$$ H_1(X^\omega,\Z)\cong \Z^{b_1(X)+r} \oplus \Z/d_1 \Z \oplus \cdots \oplus \Z/d_k \Z$$
with $1<d_1 \mid \cdots \mid d_k$,  if and only if
$$ H_1(X,\sL_\omega) \cong  \Z^r \oplus \Z/2 d_1 \Z \oplus \cdots \oplus \Z/2 d_k \Z \oplus (\Z_2)^{b_1(X)-r-k-1} .$$
\ec

Since the complement $M(\sA)$ of a complex hyperplane arrangement $\sA$ is homotopy equivalent to a minimal CW complex (see \cite{DP03} or \cite{Ran02}), this corollary gives a positive answer to the following conjecture proposed by Ishibashi, Sugawara and Yoshinaga recently. 
\begin{conj} \cite[Conjecture 3.3]{ISY22} Let $M(\sA)$ be the complement of a complex hyperplane arrangement $\sA$ in $\C^n$. For a nonzero $\omega\in H^1(M(\sA),\Z_2)$, the integral homology $H_1(M(\sA)^{\omega},\Z)$ has a non-trivial $2$-torsion if and only if the local system homology $H_1(M(\sA),\sL_{\omega})$ has a non-trivial $4$-torsion. 
\end{conj}
In particular, Corollary \ref{cor1} is compatible with the computations in \cite[Example 3.1, Example 3.2]{ISY22}.

\bc \label{cor2}  With the same assumptions and notations as in Theorem \ref{main}, the torsion part of $H_i(X, \sL_\omega)$ is either 0 or a finite direct sum of $\Z_2$ for all $i$ if and only if $H_i(X^\omega,\Z)$ is torsion-free for all $i$.
\ec

Corollary \ref{cor2}  follows directly from Theorem \ref{main} and the Universal Coefficient Theorem.
Using Corollary \ref{cor2} and Sugawara's recent work on real hyperplane arrangements \cite{Sug22}, we give a partial answer to a question asked by Yoshinaga \cite[Remark 3.8]{Yos20}.

\bc \label{application} Let $\sA$ be a real hyperplane arrangement in $\mathbb{R}^n$ and $X$ be the complexified complement of $\sA$ in $\C^n$. Fix a nonzero element $\omega\in H^1(X,\Z_2)$. 
If $\sL_\omega$ satisfies the CDO-condition (see \cite[Definition 2.1]{Sug22}), then $H_i(X^\omega,\Z)$ is torsion-free for all $i$. Moreover, $H_*(X^\omega, \Z)$ is combinatorially determined.
\ec

%Question: In \cite[Example 3.1, Example 3.2]{ISY22}, there are only 2-torsion in $H_1(X^\omega,\Z)$. Is it possible to have $p$-torsion with $p>2$, where $p$ is a prime number.

\section{Proofs}
\begin{proof}[Proof of Theorem \ref{main}]
Since the homology groups of local systems are homotopy invariants, without loss of generality, we assume that $X$ is a minimal CW complex from now on. Note that  $H_1(X,\Z)$ is torsion free in this case. Then there exists a surjective group homomorphism
  $\nu \colon \pi_1(X)\to \Z$ with the following commutative diagram 
\be \label{a}
\xymatrix{
\pi_1(X)\ar@{->>}"1,3"^{\omega}\ar@{->>}[dr]_{\nu} & & \Z_2, \\
 & \Z\ar@{->>}[ur] &
}\ee
where $\Z\to \Z_2$ is the natural quotient map.  Hence the representation of the rank one $\Z$-local system $\sL_\omega$ factors through $\nu$ and it sends the generator $1_\Z$ of $\Z$ to $-1\in \Z^\times$.

Let $X^{\nu}$ be the covering space of $X$ associated to $\nu$.
The group of covering transformations of $X^{\nu}$ is isomorphic to $\Z$ and acts on it. Consider the minimal CW-structure of $X$.  By choosing lifts of these cells of $X$ to $X^{\nu}$, we obtain a free basis for the cellular chain complex $C_{*}(X^{\nu}, \Z)$ of $ X^\nu$ as  $\Z[t^\pm]$-modules.   So $C_{*}(X^{\nu}, \Z)$ is a bounded complex of finitely generated free $\Z[t^\pm]$-modules:
$$
 \cdots  \to  C_{i+1}(X^\nu, \Z) \overset{\partial_{i}}{\to} C_i(X^\nu, \Z) \overset{\partial_{i-1}}{\to} C_{i-1}(X^\nu, \Z)  \overset{\partial_{i-2}}{\to}   \cdots \overset{\partial_0}{\to} C_0(X^\nu, \Z)  \to 0 .
$$
With the above free basis for $C_*(X^\nu,\Z)$, $\partial_i$ can be written down as a matrix, say a $(m\times n)$-matrix $\big( f_{kj}\big)_{m\times n}$ with $f_{kj}\in \Z[t^{\pm}]$.  We can choose $f_{kj}$ to be an integral valued polynomial  simultaneously. Assume that 
\begin{center}
$f_{kj} \equiv a_{kj}t+b_{kj} \mod (t^2-1)$ with $a_{kj},b_{kj}\in \Z$.
\end{center}
Since $X$ is a minimal CW complex, every entry in $\partial_i$ is divisible by $t-1$. Hence $f_{kj}$ takes the value of 0 at $t=1$, which implies $b_{kj}=-a_{kj}$.

Due to the commutative diagram (\ref{a}), $X^\nu$ is a covering space of $X^\omega$. In particular, $H_*(X^\omega,\Z)$
 can be computed by the chain complex
$$C_*(X^\nu,\Z)\otimes_{\Z[t^{\pm}]} \Z[t^{\pm}]/(t^2-1),$$
and $H_*(X,\sL_\omega)$ can be computed by the following chain complex
 $$(E_*, \alpha_*) \coloneqq C_*(X^\nu,\Z)\otimes_{\Z[t^{\pm}]} \Z[t^{\pm}]/(t+1),$$ 
both viewed as complexes of finitely generated free abelian groups. 
Note that $f_{kj}$ takes the value $-2a_{kj}$ at $t=-1$. Hence the boundary map $\alpha_i$ can be written down as a matrix  $\big(-2 a_{kj}\big)_{m\times n}$.

On the other hand, consider the direct sum $$\Z[t^{\pm}]/(t^2-1) \cong \Z\cdot 1 \oplus \Z\cdot t  $$  as $\Z$-modules. Then by this choice of basis, $\big(f_{kj}\big)_{m\times n}$ gives a $(2m\times 2n)$-matrix with entry $f_{kj}$ being replaced by $ \bigl( \begin{smallmatrix}
-a_{kj}& a_{kj}\\ a_{kj}& -a_{kj}
\end{smallmatrix} \bigr).$
By elementary row and column operations, it becomes $$  \begin{pmatrix}
(-a_{kj})_{m\times n} & 0 \\
0& 0
\end{pmatrix}_{2m\times 2n} .$$
This matrix has the same invariant divisors as of the linear map $ \frac{1}{2} \alpha_i$. It implies  
$$\tor (H_i(X^\omega,\Z))\cong \tor (H _i(E_*, \tfrac{1}{2}\alpha_*)),$$
where $\tor(-)$ denotes the torsion part of the corresponding abelian group. Since $X$ is  a minimal CW complex, $H_i(X,\Z)$ is torsion free. So the isomorphism (\ref{*}) holds for the torsion parts.

Let $\pi\colon X^\omega \to X$ be the double covering map associated to $\omega$. Then  we have $$\pi_* \C_{X^\omega}\cong \C_X\oplus (\sL_\omega\otimes \C),$$
hence $$H_i(X^\omega,\C)\cong H_i(X,\C)\oplus H_i(X,\sL_\omega\otimes \C)\cong H_i(X,\C)\oplus H_i(E_*\otimes \C,\tfrac{1}{2}\alpha_*).$$ Therefore the isomorphism (\ref{*}) also holds for the free abelian parts.
 Then the claim follows.
\end{proof}

\begin{proof}[Proof of Corollary \ref{cor1}]
Since the homology groups of local systems are homotopy invariants, without loss of generality, we assume that $X$ is a minimal CW complex.
Then $H_1(X,\Z)$ is free of rank $b_1(X)$. Theorem \ref{main} implies that $H_1(X^\omega,\Z)$ has the displayed form if and only if
$$H_1(E_*,\frac{1}{2}\alpha_*)\cong \Z^{r} \oplus \Z/d_1 \Z \oplus \cdots \oplus \Z/d_k \Z.$$
Since $X$ is connected and  $\sL_\omega$ is a non-trivial local system, $E_0\cong \Z$ and $H_0(E_*,\frac{1}{2}\alpha_*)=0$. 
Then   the invariant divisors of $\frac{1}{2}\alpha_1:E_2\to E_1$ are
$$\{ 1,\ldots,1,d_1,\ldots,d_k\},$$
where the multiplicity of $1$ is $b_1(X)-1-k-r$. One readily sees that this is equivalent to that $H_1(X,\sL_\omega)\cong H_1(E_*,\alpha_*)$ has the displayed form.
\end{proof}

\begin{proof}[Proof of Corollary \ref{application}]
To give the proof, we need to show that the torsion part of $H_i(X, \sL_\omega)$ is either 0 or a direct sum of $\Z_2$ for all $i$. The cohomology version of this claim is proved in \cite[Theorem 1.3]{Sug22}. Then one gets the homology version by the Universal Coefficient Theorem.  Corollary \ref{cor2} implies that $H_*(X^\omega,\Z)$ is torsion-free for all degrees. Hence $H_*(X^\omega, \Z)$ is combinatorially determined due to \cite[Theorem 2.4]{ISY22}.
\end{proof}

\bex[Double star arrangement] Consider the double star arrangement \cite[Figure 2]{ISY22} and we use the same notations as in \cite[Example 3.2]{ISY22}.
Take $$\omega =e_8+e_9+e_{10}+e_{12}+e_{13}+e_{14}+e_{15}.$$ It is easy to check that the CDO-condition holds for $\sL_\omega$ on the hyperplane $H_6$.
Then we have $H_1(M(\sA_{\mathcal{DS}}),\sL_\omega)\cong  \Z_2^9$ (see \cite[Theorem 1.3]{Sug22}), hence 
$$
H_1(M(\sA_{\mathcal{DS}})^\omega,\Z)\cong \Z^{10}.
$$ 
\eex

\section*{Acknowledgments} We thank Masahiko Yoshinaga for valuable comments. We also thank the referee for useful comments and suggestions. The first named author is supported by NSFC grant No. 11901467. The second named author is  supported by National Key Research and Development Project SQ2020YFA070080, the starting grant from University of Science and Technology of China, NSFC grant No. 12001511, the Project of Stable Support for Youth Team in Basic Research Field, CAS (YSBR-001),  the project ``Analysis and Geometry on Bundles" of Ministry of Science and Technology of the People's Republic of China and  Fundamental Research Funds for the Central Universities.

\bibliographystyle{amsalpha}

\end{document}